\documentclass[a4paper,11pt,fleqn]{article}
\usepackage[obeyspaces,hyphens,spaces]{url}
\usepackage[noadjust]{cite}
\usepackage[dvipsnames]{xcolor}
\usepackage{amsmath}
\usepackage{amssymb}
\usepackage{amsfonts}
\usepackage{charter}
\usepackage{oplotsymbl}
\usepackage{etoolbox}
\usepackage{mathtools}
\usepackage{euscript}
\usepackage{mathrsfs} %pour mathscr
\usepackage{pifont} %for dingolist
\usepackage{srcltx} %for forward references in xdvi
\usepackage{multirow}
\allowdisplaybreaks %allows align to split
\usepackage{tikz,tkz-euclide,pgfplots}
\usetikzlibrary{arrows,decorations}
\usepackage[usenames,dvipsnames]{pstricks}
\usepackage{pstricks-add}
\usepackage{pst-grad}
\usepackage{pst-plot} 
\usepackage{pst-eucl} 
\psset{algebraic}
\usepackage{theorem} 
\usepackage{euscript}
\usepackage{exscale,relsize}
\usepackage{epic,eepic}
\usepackage{pgfplots}
\usepackage{makeidx}
\usepackage{multicol}
\usepackage{graphicx}
\definecolor{labelkey}{HTML}{0455BF}
\definecolor{refkey}{rgb}{0,0.6,0.0}
\definecolor{lblue}{HTML}{0455BF}
\definecolor{dgreen}{HTML}{02724A}
\definecolor{myellow}{HTML}{D97904}
\definecolor{orng}{HTML}{D35400}
\definecolor{dred}{HTML}{D90404}
\definecolor{pblue}{rgb}{0.1176,0.5647,1}
\definecolor{pgreen}{rgb}{0.1961,0.8039,0.1961}
\definecolor{pred}{rgb}{1.0,0.2706,0.0}
\definecolor{fred}{rgb}{0.93,0.51,0.93}
\definecolor{pyellow}{rgb}{1.0,0.6471,0.0}
\usepackage{upref}%upright number when using ref in theorems
\usepackage{hyperref}
\hypersetup{colorlinks=true,linktocpage=true,linkcolor=lblue,%
citecolor=dgreen,urlcolor=dred}

\renewcommand\familydefault{\rmdefault}
\DeclareMathAlphabet{\mathrm}{OT1}{\familydefault}{m}{n}
\makeatletter
\def\operator@font{\mathgroup\symoperators\rm}
\makeatother
\renewcommand{\leq}{\ensuremath{\leqslant}}

\newcommand{\minimize}[2]{\underset{#1}{\operatorname{minimize}}
\;\:#2}

\newcommand{\scal}[2]{{\langle{{#1}\mid{#2}}\rangle}}
\newcommand{\sscal}[2]{{\big\langle{{#1}\mid{#2}}\big\rangle}}

\newcommand{\menge}[2]{\big\{{#1}~|~{#2}\big\}} 
 
\newcommand{\GGG}{\ensuremath{\boldsymbol{\mathcal{G}}}}
\newcommand{\HHH}{\ensuremath{\boldsymbol{\mathcal{H}}}}
\newcommand{\HH}{\ensuremath{{\mathcal{H}}}}

\newcommand{\GG}{\ensuremath{{\mathcal{G}}}}

\newcommand{\RR}{\ensuremath{\mathbb{R}}}
\newcommand{\NN}{\ensuremath{\mathbb{N}}}

\newcommand{\pinf}{\ensuremath{{{+}\infty}}}

\newcommand{\RX}{\ensuremath{\left]{-}\infty,{+}\infty\right]}}

\newcommand{\RPP}{\ensuremath{\left]0,{+}\infty\right[}}

\newcommand{\emp}{\ensuremath{\varnothing}}
\newcommand{\Sum}{\ensuremath{\displaystyle\sum}}

\newcommand{\ID}{\ensuremath{\mathbf{Id}}}

\newcommand{\prox}{\mathrm{prox}}
\newcommand{\proj}{\mathrm{proj}}
\newcommand{\prob}{\mathsf{Prob}}

\usepackage{theorem}
\theoremstyle{plain}{\theorembodyfont{\rmfamily}%
\newtheorem{algorithm}{Algorithm}}

\newtheorem{theorem}{Theorem}

\theoremstyle{plain}{\theorembodyfont{\rmfamily}%
}
\theoremstyle{plain}{\theorembodyfont{\rmfamily}%
}
\theoremstyle{plain}{\theorembodyfont{\rmfamily}%
}
\theoremstyle{plain}{\theorembodyfont{\rmfamily}%
\newtheorem{problem}[theorem]{Problem}}
\theoremstyle{plain}{\theorembodyfont{\rmfamily}%
}
\usepackage{enumitem}
\setlist[enumerate]{itemsep=2pt}
\setlist[itemize]{itemsep=2pt}

\newcommand*\mute{{\mkern 2mu\cdot\mkern 2mu}}
\usepackage{authblk}

\newcommand{\email}[1]{\href{mailto:#1}{\nolinkurl{#1}}}
\begin{document}
\title{Block-Activated Algorithms for Multicomponent Fully
Nonsmooth Minimization\thanks{M. N. B\`ui and P. L. Combettes were
supported by the National Science Foundation under grant
CCF-1715671 and Z. C. Woodstock by the National Science Foundation
under grant DGE-1746939.}}
\author{Minh N. B\`ui, Patrick L. Combettes, and 
Zev C. Woodstock}
\affil{North Carolina State University,
Department of Mathematics, Raleigh, NC 27695-8205, USA\\
\email{mnbui@ncsu.edu}, \email{plc@math.ncsu.edu}, 
\email{zwoodst@ncsu.edu}}

\date{~}
\maketitle

\begin{abstract}
Under consideration are multicomponent minimization problems
involving a separable nonsmooth convex function penalizing the
components individually, and nonsmooth convex coupling terms
penalizing linear mixtures of the components. We investigate
block-activated proximal algorithms for solving such problems,
i.e., algorithms which, at each iteration, need to use only a block
of the underlying functions, as opposed to all of them as in
standard methods. For smooth coupling functions, several
block-activated algorithms exist and they are well understood. By
contrast, in the fully nonsmooth case, few block-activated methods
are available and little effort has been devoted to assessing them.
Our goal is to shed more light on the
implementation, the features, and the behavior of these algorithms,
compare their merits, and provide machine learning and image
recovery experiments illustrating their performance.
\end{abstract}

\section{Introduction} 
\label{sec:1}
The goal of many signal processing and machine learning tasks is to
exploit the observed data and the prior knowledge to produce a
solution that represents information of interest. In
this process of extracting information from data, structured convex
optimization has established itself as an effective modeling and
algorithmic framework; see, for instance,
\cite{Bach12,Boyd10,Cham16,Banf11,Glow16}. In
state-of-the-art applications, the sought solution is often
a tuple of vectors which reside in different spaces
\cite{Argy12,Aujo05,Biou12,Bric19,Nume19,Sign21,Ejst20,%
Darb20,Hint06}. The following multicomponent minimization
formulation captures such problems. It consists of a separable term
penalizing the components individually, and of coupling terms
penalizing linear mixtures of the components.

\begin{problem}
\label{prob:1}
Let $(\HH_i)_{1\leq i\leq m}$ and $(\GG_k)_{1\leq k\leq p}$ be 
Euclidean spaces. For every $i\in\{1,\ldots,m\}$ and every
$k\in\{1,\ldots,p\}$, let $f_i\colon\HH_i\to\RX$ and
$g_k\colon\GG_k\to\RX$ be proper lower semicontinuous convex
functions, and let $L_{k,i}\colon\HH_i\to\GG_k$ be a linear 
operator. The objective is to
\begin{equation}
\label{e:prob1}
\minimize{x_1\in\HH_1,\ldots,x_m\in\HH_m}
{\underbrace{\sum_{i=1}^mf_i(x_i)}_{\text{separable term}}
+\sum_{k=1}^p\underbrace{g_k\Bigg(\sum_{i=1}^mL_{k,i}
x_i\Bigg)}_{\text{$k$th coupling term}}}.
\end{equation}
\end{problem}

To solve Problem~\ref{prob:1} reliably without adding restrictions
(for instance, smoothness or strong convexity of
some functions involved in the model), we focus on flexible
proximal algorithms that have the following features:
\begin{dingautolist}{192}
\item
\label{li:1}
{\bfseries Nondifferentiability:}
None of the functions $f_1,\ldots,$ $f_m,g_1,\ldots,g_p$ 
needs to be differentiable.
\item
\label{li:2}
{\bfseries Splitting:} The functions
$f_1,\ldots,f_m,g_1,\ldots,g_p$ and the
linear operators are activated separately.
\item
\label{li:3}
{\bfseries Block activation:}
Only a block of the functions $f_1,\ldots,$ $f_m,g_1,\ldots,g_p$ is
activated at each iteration. This is in contrast with most
splitting methods which require full activation, i.e., that all the
functions be used at every iteration.
\item
\label{li:5}
{\bfseries Operator norms:} Bounds on the norms of the linear
operators involved in Problem~\ref{prob:1} are not assumed since
they can be hard to compute. 
\item
\label{li:4}
{\bfseries Convergence guarantee:}
The algorithm produces a sequence which converges (possibly 
almost surely) to a solution to Problem~\ref{prob:1}.
\end{dingautolist}

In view of features \ref{li:1} and \ref{li:2}, the
algorithms of interest should activate the functions
$f_1,\ldots,f_m,g_1,\ldots,g_p$ via their proximity
operators (even if some functions happened
to be smooth, proximal activation is often preferable
\cite{Bric19,Siim19}). The motivation for \ref{li:2} is that
proximity 
operators of composite functions are typically not known
explicitly. Feature~\ref{li:3} is geared towards
current large-scale problems. In such scenarios, memory and
computing power limitations make the execution of standard
proximal splitting algorithms, which require activating all
the functions at each iteration, inefficient or
simply impossible. We must therefore turn our attention to
algorithms which employ only
blocks of functions $(f_i)_{i\in I_n}$ and $(g_k)_{k\in K_n}$ at 
iteration $n$. If
the functions $(g_k)_{1\leq k\leq p}$ were all smooth, one could
use block-activated versions of the forward-backward algorithm
proposed in \cite{Siop15,Salz20} and the references therein; in
particular, when $m=1$, methods such as those of
\cite{Anon21,Defa14,Mish20,Bach17} would be pertinent. 
As noted in \cite[Remark~5.10(iv)]{Siop15}, another candidate of
interest could be the randomly block-activated algorithm of
\cite[Section~5.2]{Siop15}, which leads to block-activated versions
of several primal-dual methods (see \cite{Pesq15} for detailed
developments and \cite{Cham18} for an inertial version when $m=1$).
However, this approach violates \ref{li:5} as it
imposes bounds on the proximal scaling parameters which depend on
the norms of the linear operators. Finally, \ref{li:4}
rules out methods that guarantee merely minimizing sequences or
ergodic convergence.

To the best of our knowledge, there are two primary
methods that fulfill \ref{li:1}--\ref{li:4}:
\begin{itemize}
\item
Algorithm~\ref{algo:A}:
The stochastic primal-dual Douglas--Rachford algorithm
of \cite{Siop15}.
\item
Algorithm~\ref{algo:B}:
The deterministic primal-dual projective splitting algorithm of
\cite{MaPr18}. 
\end{itemize}
In the case of smooth coupling functions $(g_k)_{1\leq k\leq p}$,
in \eqref{e:prob1}, extensive numerical experience has been
accumulated to understand the behavior of block-activated methods,
especially in the case of stochastic gradient methods. By contrast,
to date, very few numerical experiments with the recent, fully
nonsmooth Algorithms~\ref{algo:A} and \ref{algo:B} have been
conducted and no comparison of their merits and performance has
been undertaken. Thus far, Algorithm~\ref{algo:A} has been employed
only in the context of machine learning (see also the variant of
\ref{algo:A} in \cite{Bric19} for partially smooth problems). On
the other
hand, Algorithm~\ref{algo:B} has been used in image recovery in
\cite{Siim19}, but only in full activation mode, and in 
feature selection in \cite{John20}, but with $m=1$.

{\bfseries Contributions and novelty:}
This paper investigates for the first time the use of
block-activated methods in fully nonsmooth multivariate
minimization problems. It sheds more light on the implementation,
the features, and the behavior of 
Algorithms~\ref{algo:A} and \ref{algo:B}, compares their merits,
and provides experiments illustrating their performance. 

{\bfseries Outline:}
Algorithms~\ref{algo:A} and \ref{algo:B} are presented in
Section~\ref{sec:2}. In Section~\ref{sec:3}, we analyze and compare
their features, implementation, and asymptotic properties. This
investigation is complemented in Section~\ref{sec:4} by numerical
experiments in the context of machine learning and image recovery. 

\section{Block-Activated Algorithms for Problem~\ref{prob:1}}
\label{sec:2}

The subdifferential, the conjugate, and the proximity operator
of a proper lower semicontinuous convex function $f\colon\HH\to\RX$
are denoted by $\partial f$, $f^*$, and $\prox_f$, respectively.
Let us consider the setting of Problem~\ref{prob:1} and let us set
$\HHH=\HH_1\times\cdots\times\HH_m$ and
$\GGG=\GG_1\times\cdots\times\GG_p$. A generic element
in $\HHH$ is denoted by $\boldsymbol{x}=(x_i)_{1\leq i\leq m}$ and
a generic element in $\GGG$ by 
$\boldsymbol{y}=(y_k)_{1\leq k\leq p}$.

As discussed in Section~\ref{sec:1}, two primary algorithms 
fulfill requirements~\ref{li:1}--\ref{li:4}. Both 
operate in the product space $\HHH\times\GGG$. The first one
employs random activation of the blocks. To present it, let us
introduce 
\begin{equation}
\label{e:LL}
\begin{cases}
\boldsymbol{L}\colon\HHH\to\GGG\colon
\boldsymbol{x}\mapsto\bigg(\Sum_{i=1}^m{L}_{1,i}x_i,
\ldots,\Sum_{i=1}^m{L}_{p,i}x_i\bigg)\\
\boldsymbol{V}=
\menge{(\boldsymbol{x},\boldsymbol{y})\in\HHH\times\GGG}
{\boldsymbol{y}=\boldsymbol{L}\boldsymbol{x}}\\
\boldsymbol{F}\colon\HHH\times\GGG\to\RX\\
~~~~\hspace{1.9mm}(\boldsymbol{x},\boldsymbol{y})
\mapsto\sum_{i=1}^mf_i(x_i)+\sum_{k=1}^pg_k(y_k).
\end{cases}
\end{equation}
Then \eqref{e:prob1} is equivalent to 
\begin{equation}
\label{e:094}
\minimize{(\boldsymbol{x},\boldsymbol{y})\in
\boldsymbol{V}}{\boldsymbol{F}(\boldsymbol{x},\boldsymbol{y})}.
\end{equation}
The idea is then to apply the Douglas--Rachford algorithm in block
form to this problem \cite{Siop15}. To this end, we need 
$\prox_{\boldsymbol{F}}$ and 
$\prox_{\iota_{\boldsymbol{V}}}=\proj_{\boldsymbol{V}}$. Note that
$\prox_{\boldsymbol{F}}\colon(\boldsymbol{x},\boldsymbol{y})
\mapsto((\prox_{f_i}x_i)_{1\leq i\leq m},
(\prox_{g_k}y_k)_{1\leq k\leq p})$.
Now let $\boldsymbol{x}\in\HHH$ and
$\boldsymbol{y}\in\GGG$, and set
$\boldsymbol{t}=(\ID+\boldsymbol{L}^*
\boldsymbol{L})^{-1}(\boldsymbol{x}+
\boldsymbol{L}^*\boldsymbol{y})$ and 
$\boldsymbol{s}=(\ID+\boldsymbol{L}\boldsymbol{L}^*)^{-1}
(\boldsymbol{L}\boldsymbol{x}-\boldsymbol{y})$.
Then 
\begin{equation}
\label{e:o7}
\proj_{\boldsymbol{V}}(\boldsymbol{x},\boldsymbol{y})
=(\boldsymbol{t},\boldsymbol{L}\boldsymbol{t})=
(\boldsymbol{x}-\boldsymbol{L}^*\boldsymbol{s},
\boldsymbol{y}+\boldsymbol{s}),
\end{equation}
and we write it coordinate-wise as
\begin{equation}
\label{e:o8}
\proj_{\boldsymbol{V}}(\boldsymbol{x},\boldsymbol{y})=
\big(Q_1(\boldsymbol{x},\boldsymbol{y}),\ldots,Q_{m+p}
(\boldsymbol{x},\boldsymbol{y})\big).
\end{equation}
Thus, given $\gamma\in\RPP$, $\boldsymbol{z}_0\in\HHH$, and
$\boldsymbol{y}_0\in\GGG$, the standard Douglas--Rachford 
algorithm for \eqref{e:094} is 
\begin{equation}
\begin{array}{l}
\text{for}\;n=0,1,\ldots\\
\left\lfloor
\begin{array}{l}
\lambda_n\in\left]0,2\right[\\
\text{for every}\;i\in\{1,\ldots,m\}\\
\left\lfloor
\begin{array}{l}
x_{i,n+1}=Q_i(\boldsymbol{z}_n,\boldsymbol{y}_n)\\
z_{i,n+1}=z_{i,n}+\lambda_n
\big(\prox_{\gamma f_i}(2x_{i,n+1}-z_{i,n})-x_{i,n+1}\big)
\end{array}
\right.
\\
\text{for every}\;k\in\{1,\ldots,p\}\\
\left\lfloor
\begin{array}{l}
w_{k,n+1}=Q_{m+k}(\boldsymbol{z}_n,\boldsymbol{y}_n)\\
y_{k,n+1}=y_{k,n}+\lambda_n
\big(\prox_{\gamma g_k}(2w_{k,n+1}-y_{k,n})-w_{k,n+1}\big).
\end{array}
\right.\\[6mm]
\end{array}
\right.
\end{array}
\end{equation}
The block-activated version of this algorithm is as follows.

\begin{algorithm}[\cite{Siop15}]
\label{algo:A}
Let $\gamma\in\RPP$,
let $\boldsymbol{x}_0$ and $\boldsymbol{z}_0$ be $\HHH$-valued
random variables (r.v.), let $\boldsymbol{y}_0$ and
$\boldsymbol{w}_0$ be $\GGG$-valued r.v. Iterate
\[
\hskip -3mm
\begin{array}{l}
\text{for}\;j=1,\ldots,m+p\\
\left\lfloor
\text{compute}\;Q_j\;\text{as in \eqref{e:o7}--\eqref{e:o8}}
\right.
\\
\text{for}\;n=0,1,\ldots\\
\left\lfloor
\begin{array}{l}
\lambda_n\in\left]0,2\right[\\
\text{select randomly}\;\emp\neq I_n\subset\{1,\ldots,m\}
\;\text{and}\;\emp\neq K_n\subset\{1,\ldots,p\}\\
\text{for every}\;i\in I_n\\
\left\lfloor
\begin{array}{l}
x_{i,n+1}=Q_i(\boldsymbol{z}_n,\boldsymbol{y}_n)\\
z_{i,n+1}=z_{i,n}+
\lambda_n\big(\prox_{\gamma f_i}(2x_{i,n+1}-z_{i,n})
-x_{i,n+1}\big)
\end{array}
\right.
\\
\text{for every}\;i\in \{1,\ldots,m\}\smallsetminus I_n\\
\left\lfloor
\begin{array}{l}
(x_{i,n+1},z_{i,n+1})=(x_{i,n},z_{i,n})
\end{array}
\right.
\\
\text{for every}\;k\in K_n\\
\left\lfloor
\begin{array}{l}
w_{k,n+1}=Q_{m+k}(\boldsymbol{z}_n,\boldsymbol{y}_n)\\
y_{k,n+1}=y_{k,n}+\lambda_n
\big(\prox_{\gamma g_k}(2w_{k,n+1}-y_{k,n})-w_{k,n+1}\big)
\end{array}
\right.
\\
\text{for every}\;k\in \{1,\ldots,p\}\smallsetminus K_n\\
\left\lfloor
\begin{array}{l}
(w_{k,n+1},y_{k,n+1})=(w_{k,n},y_{k,n}).
\end{array}
\right.\\[2mm]
\end{array}
\right.
\end{array}
\]
\end{algorithm}

The second algorithm operates by projecting onto hyperplanes which
separate the current iterate from the set 
$\boldsymbol{\mathsf{Z}}$ of Kuhn--Tucker points of
Problem~\ref{prob:1}, i.e., the points 
$\widetilde{\boldsymbol{x}}\in\HHH$ and
$\widetilde{\boldsymbol{v}}^*\in\GGG$ such that 
\begin{equation}
\begin{cases}
\label{e:Z}
(\forall i\in\{1,\ldots,m\})\quad{-}
\sum_{k=1}^pL_{k,i}^*\widetilde{v}_k^*\in
\partial f_i(\widetilde{x}_i)\\
(\forall k\in\{1,\ldots,p\})\quad
\sum_{i=1}^mL_{k,i}\widetilde{x}_i\in\partial g_k^*
(\widetilde{v}_k^*).
\end{cases}
\end{equation}
This process is explained in Fig.~\ref{fig:8}.

\begin{algorithm}[\cite{MaPr18}]
\label{algo:B}
Set $I_0=\{1,\ldots,m\}$ and $K_0=\{1,\ldots,p\}$. For every 
$i\in I_0$ and every $k\in K_0$, let
$\{\gamma_i,\mu_k\}\subset\RPP$, $x_{i,0}\in\HH_i$, and
$v_{k,0}^*\in\GG_k$. Iterate
\[
\label{e:a1}
\begin{array}{l}
\text{for}\;n=0,1,\ldots\\
\left\lfloor
\begin{array}{l}
\lambda_n\in\left]0,2\right[\\
\text{if}\;n>0\\
\left\lfloor
\begin{array}{l}
\text{select}\;\emp\neq I_n\subset I_0
\;\text{and}\;\emp\neq K_n\subset K_0
\end{array}
\right.
\\
\text{for every}\;i\in I_n\\
\left\lfloor
\begin{array}{l}
x^*_{i,n}=x_{i,n}-\gamma_i\sum_{k=1}^pL_{k,i}^*v_{k,n}^*\\
a_{i,n}=\prox_{\gamma_if_i}x_{i,n}^*\\
a_{i,n}^*=\gamma_i^{-1}(x_{i,n}^*-a_{i,n})
\end{array}
\right.\\
\text{for every}\;i\in I_0\smallsetminus I_n\\
\left\lfloor
\begin{array}{l}
(a_{i,n},a_{i,n}^*)=(a_{i,n-1},a_{i,n-1}^*)
\end{array}
\right.\\
\text{for every}\;k\in K_n\\
\left\lfloor
\begin{array}{l}
y_{k,n}^*=\mu_kv_{k,n}^*+\sum_{i=1}^mL_{k,i}x_{i,n}\\
b_{k,n}=\prox_{\mu_kg_k}y_{k,n}^*\\
b^*_{k,n}=\mu_k^{-1}(y_{k,n}^*-b_{k,n})\\
t_{k,n}=b_{k,n}-\sum_{i=1}^mL_{k,i}a_{i,n}
\end{array}
\right.\\
\text{for every}\;k\in K_0\smallsetminus K_n\\
\left\lfloor
\begin{array}{l}
(b_{k,n},b^*_{k,n})=(b_{k,n-1},b^*_{k,n-1})\\
t_{k,n}=b_{k,n}-\sum_{i=1}^mL_{k,i}a_{i,n}
\end{array}
\right.\\
\text{for every}\;i\in I_0\\
\left\lfloor
\begin{array}{l}
t^*_{i,n}=a^*_{i,n}+\sum_{k=1}^pL_{k,i}^*b^*_{k,n}
\end{array}
\right.\\
\tau_n=\sum_{i=1}^m\|t_{i,n}^*\|^2+\sum_{k=1}^p\|t_{k,n}\|^2\\
\text{if}\;\tau_n>0\\
\left\lfloor
\begin{array}{l}
\pi_n=\textstyle\sum_{i=1}^m\big(\scal{x_{i,n}}{t^*_{i,n}}-
\scal{a_{i,n}}{a^*_{i,n}}\big)\\
\textstyle\quad\;+\sum_{k=1}^p\big(\sscal{t_{k,n}}{v_{k,n}^*}
-\sscal{b_{k,n}}{b^*_{k,n}}\big)
\end{array}
\right.\\
\text{if}\;\tau_n>0\;\text{and}\;\pi_n>0\\
\left\lfloor
\begin{array}{l}
\theta_n=\lambda_n\pi_n/\tau_n\\
\text{for every}\;i\in I_0\\
\left\lfloor
\begin{array}{l}
x_{i,n+1}=x_{i,n}-\theta_nt^*_{i,n}
\end{array}
\right.\\
\text{for every}\;k\in K_0\\
\left\lfloor
\begin{array}{l}
v_{k,n+1}^*=v_{k,n}^*-\theta_nt_{k,n}
\end{array}
\right.\\[1mm]
\end{array}
\right.
\\
\text{else}\\
\left\lfloor
\begin{array}{l}
\text{for every}\;i\in I_0\\
\left\lfloor
\begin{array}{l}
x_{i,n+1}=x_{i,n}
\end{array}
\right.\\
\text{for every}\;k\in K_0\\
\left\lfloor
\begin{array}{l}
v_{k,n+1}^*=v_{k,n}^*.
\end{array}
\right.\\[1mm]
\end{array}
\right.\\[8mm]
\end{array}
\right.
\end{array}
\]
\end{algorithm}

\begin{figure}[ht]
\centering
\begin{tikzpicture}[scale=1.0]
\tkzInit[xmin=-1.2,xmax=8,ymin=-1,ymax=5.1]
\tkzDrawX[line width=0.06cm,%
noticks,%
label={\normalsize${\HHH}$},%
-latex']
\tkzDrawY[line width=0.06cm,%
noticks,%
label={\normalsize${\GGG}$},%
-latex']
%ellipse
\begin{scope}[rotate around={25:(0.58,1.6)}]
\draw[color=black,fill opacity=0.3,
line width=0.06cm,fill=lblue] (3.18,1.54)
ellipse (2.45cm and 1.42cm);
\end{scope}
%lines
\draw[line width=0.04cm,-latex'] (6.8,0.7) -- (5.7,1.7);
\draw[line width=0.06cm,color=fred] (0,0.95) -- (0,4.35);
\draw[line width=0.06cm,color=pblue] (0.6,0) -- (5.25,0);
\draw[line width=0.06cm,color=dgreen] (3.1,-1.2) -- (7.9,4.1);
% nodes
\node at (6.9,3.8)
{\color{dgreen}\normalsize${\boldsymbol{\mathsf{H}}_n}$};
\node at (3.5,3.4) {\normalsize${\boldsymbol{\mathsf{Z}}}$};
\tkzDefPoint(6.8,0.0){xn}
\tkzDefPoint(5.7,0.0){xn1}
\tkzDefPoint(0.0,0.7){vn}
\tkzDefPoint(0.0,1.7){vn1}
\tkzDefPoint(6.8,0.7){Xn}
\tkzDefPoint(5.7,1.7){Xn1}
\tkzDefPoint(3.0,0){P}
\tkzDefPoint(0,3.0){D}
%\tkzLabelPoint[below](xn){$\boldsymbol{x}_n$}
%\tkzLabelPoint[below](xn1){$\boldsymbol{x}_{n+1}$}
%\tkzLabelPoint[left](vn){$\boldsymbol{v}_n^*$}
%\tkzLabelPoint[left](vn1){$\boldsymbol{v}_{n+1}^*$}
\tkzLabelPoint[below](P){\color{pblue}\normalsize$\mathscr{P}$}
\tkzLabelPoint[left](D){\color{fred}\normalsize$\mathscr{D}$}
\tkzDrawPoints(xn,xn1,vn,vn1,Xn,Xn1)
%\tkzDrawPoints(Xn,Xn1)
\tikzset{arrow coord style/.style={dashed,
line width=0.02cm,
color=black}}
\tikzset{xcoord style/.style={
font=\normalsize,text height=1ex,
inner sep = 0pt,
outer sep = 0pt,
below=3pt}}
\tikzset{ycoord style/.style={
font=\normalsize,text height=1ex,
inner sep = 0pt,
outer sep = 0pt,
left=4pt}}
\tkzShowPointCoord[xlabel={$\boldsymbol{x}_n$},
ylabel=$\boldsymbol{v}_n^*$](Xn)
\tkzShowPointCoord[xlabel={$\boldsymbol{x}_{n+1}$},
ylabel=$\boldsymbol{v}_{n+1}^*$](Xn1)
\end{tikzpicture}
\caption{Let $\mathscr{P}$ be the set of solutions to
Problem~\ref{prob:1} and let $\mathscr{D}$ be the set of 
solutions to its dual. Then the Kuhn--Tucker set
$\boldsymbol{\mathsf{Z}}$ is a subset of 
$\mathscr{P}\times\mathscr{D}$. At iteration $n$, the 
proximity operators of blocks of functions $(f_i)_{i\in I_n}$ and
$(g_k)_{k\in K_n}$ are used to construct a hyperplane 
$\boldsymbol{\mathsf{H}}_n$ separating
the current primal-dual iterate 
$(\boldsymbol{x}_n,\boldsymbol{v}_n^*)$ from
$\boldsymbol{\mathsf{Z}}$, and the update
$(\boldsymbol{x}_{n+1},\boldsymbol{v}_{n+1}^*)$ is obtained as its
projection onto $\boldsymbol{\mathsf{H}}_n$ \cite{MaPr18}.}
\label{fig:8}
\end{figure}

\section{Asymptotic Behavior and Comparisons}
\label{sec:3}

Let us first state the convergence results available for
Algorithms~\ref{algo:A} and \ref{algo:B}. We make the standing
assumption that $\boldsymbol{\mathsf{Z}}\neq\emp$ (see
\eqref{e:Z}), which implies that the solution set $\mathscr{P}$ of
Problem~\ref{prob:1} is nonempty.

\begin{theorem}[\cite{Siop15}]
\label{t:2}
In the setting of Algorithm~\ref{algo:A}, define, for every
$n\in\NN$ and every $j\in\{1,\ldots,m+p\}$,
\begin{equation}
\label{e:eps}
\varepsilon_{j,n}=
\begin{cases}
1,&\text{if}\;j\in I_n\;\text{or}\;j-m\in K_n;\\
0,&\text{otherwise}.
\end{cases}
\end{equation}
Suppose that the following hold:
\begin{enumerate}
\item
$\inf_{n\in\NN}\lambda_n>0$ and $\sup_{n\in\NN}\lambda_n<2$.
\item
The r.v. $(\boldsymbol{\varepsilon}_n)_{n\in\NN}$ are 
identically distributed.
\item
For every $n\in\NN$, the r.v. $\boldsymbol{\varepsilon}_n$ and
$(\boldsymbol{z}_j,\boldsymbol{y}_j)_{0\leq j\leq n}$ are mutually
independent.
\item
$(\forall j\in\{1,\ldots,m+p\})$
$\prob[\varepsilon_{j,0}=1]>0$.
\end{enumerate}
Then $(\boldsymbol{x}_n)_{n\in\NN}$ converges almost surely
to a $\mathscr{P}$-valued r.v. 
\end{theorem}

\begin{theorem}[\cite{MaPr18}]
\label{t:1}
In the setting of Algorithm~\ref{algo:B}, suppose that the
following hold: 
\begin{enumerate}
\item
\label{t:1a}
$\inf_{n\in\NN}\lambda_n>0$ and
$\sup_{n\in\NN}\lambda_n<2$.
\item
\label{t:1b}
There exists $T\in\NN$ such that, for every $n\in\NN$,
$\bigcup_{j=n}^{n+T}I_j=\{1,\ldots,m\}$ and
$\bigcup_{j=n}^{n+T}K_j=\{1,\ldots,p\}$.
\end{enumerate}
Then $(\boldsymbol{x}_n)_{n\in\NN}$ converges to a point in 
$\mathscr{P}$.
\end{theorem}

Let us compare Algorithms~\ref{algo:A} and \ref{algo:B}.
\begin{enumerate}[label={\rm \alph*/}]
\item
\label{r:ABi}
{\bfseries Auxiliary tasks:}
\ref{algo:A} requires the construction and storage
of the operators $(Q_j)_{1\leq j\leq m+p}$ of
\eqref{e:o7}--\eqref{e:o8},
which can be quite demanding as they
involve inversion of a linear operator acting on the product space
$\HHH$ or $\GGG$. By contrast, \ref{algo:B} does not
require such tasks.
\item
\label{r:ABii}
{\bfseries Proximity operators:}
Both algorithms are block-activated: only the blocks of functions 
$(f_i)_{i\in I_n}$ and $(g_k)_{k\in K_n}$ need to be activated at
iteration $n$.
\item
\label{r:ABiii}
{\bfseries Linear operators:}
In \ref{algo:A}, the operators $(Q_i)_{i\in I_n}$ and
$(Q_{m+k})_{k\in K_n}$ selected at iteration $n$ are evaluated at
$(z_{1,n},\ldots,z_{m,n},y_{1,n},\ldots,y_{p,n})\in
\HHH\times\GGG$. On the other hand, \ref{algo:B}
activates the local operators
$L_{k,i}\colon\HH_i\to\GG_k$ and $L^*_{k,i}\colon\GG_k\to\HH_i$
once or twice, depending on whether they are selected. 
For instance, if we set $N=\dim\HHH$ and $M=\dim\GGG$ and
if all the linear operators are implemented in matrix form,
then the corresponding load per iteration in full activation mode
of \ref{algo:A} is $\mathcal{O}((M+N)^2)$ versus
$\mathcal{O}(MN)$ in \ref{algo:B}.
\item
\label{r:ABiv}
{\bfseries Activation scheme:}
As \ref{algo:A} selects the blocks randomly,
the user does not have complete control of the computational load
of an iteration, whereas the load of \ref{algo:B} is more
predictable because of its deterministic activation scheme.
\item
\label{r:ABv}
{\bfseries Parameters:}
A single scale parameter $\gamma$ is used in
\ref{algo:A}, while \ref{algo:B}
allows the proximity operators to have their own scale parameters
$(\gamma_1,\ldots,\gamma_m,\mu_1,\ldots,\mu_p)$.
This gives \ref{algo:B} more flexibility, but more effort
may be needed \emph{a priori} to find efficient parameters.
Further, in both algorithms, there is no restriction on the
parameter values. 
\item
\label{r:ABvi}
{\bfseries Convergence:}
\ref{algo:B} guarantees sure convergence under the mild
sweeping condition \ref{t:1b} in Theorem~\ref{t:1}, while
\ref{algo:A} guarantees only almost sure convergence.
\item
\label{r:ABvii}
{\bfseries Other features:}
Although this point is omitted for brevity, unlike
\ref{algo:A}, \ref{algo:B} can be executed
asynchronously with iteration-dependent scale parameters
\cite{MaPr18}.
\end{enumerate}

\begin{figure}%[b]
\centering
\begin{tikzpicture}[scale=1.00]
\begin{axis}[height=7.2cm,width=13.9cm, legend columns=2 
cell align={left}, xmin=0, xmax=5000, ymin=-35, ymax=0.00]
\addplot [ultra thick,mark=none, densely dotted,color=pgreen]
table[x={norm-iters}, y={DR-err}] {figures2/ex1/activation10.txt};
\addplot [very thick, mark=none, color=pgreen]
table[x={norm-iters}, y={CE-err}] {figures2/ex1/activation10.txt};
\addplot [ultra thick,mark=none,dotted,color=pblue]
table[x={norm-iters}, y={DR-err}] {figures2/ex1/activation7.txt};
\addplot [very thick, mark=none, color=pblue]
table[x={norm-iters}, y={CE-err}] {figures2/ex1/activation7.txt};
\addplot [thick, mark=none,densely dashed, color=pyellow]
table[x={norm-iters}, y={DR-err}] {figures2/ex1/activation4.txt};
\addplot [very thick, mark=none, color=pyellow]
table[x={norm-iters}, y={CE-err}] {figures2/ex1/activation4.txt};
\addplot [thick, mark=none,dashed, color=pred]
table[x={norm-iters}, y={DR-err}] {figures2/ex1/activation1.txt};
\addplot [very thick, mark=none, color=pred]
table[x={norm-iters}, y={CE-err}] {figures2/ex1/activation1.txt};
\legend{%
Alg.~\ref{algo:A}-1.0,
Alg.~\ref{algo:B}-1.0,
Alg.~\ref{algo:A}-0.7,
Alg.~\ref{algo:B}-0.7,
Alg.~\ref{algo:A}-0.4,
Alg.~\ref{algo:B}-0.4,
Alg.~\ref{algo:A}-0.1,
Alg.~\ref{algo:B}-0.1}
\end{axis}
\end{tikzpicture}
\caption{Normalized error
$20\log_{10}(\|\boldsymbol{x}_n-\boldsymbol{x}_{\infty}\|
/\|\boldsymbol{x}_0-\boldsymbol{x}_{\infty}\|)$ (dB), averaged
over $20$ runs, versus epoch count in Experiment 1. The variations
around the averages were not significant. The computational load
per epoch for both algorithms is comparable.}
\label{fig:1}
\end{figure}
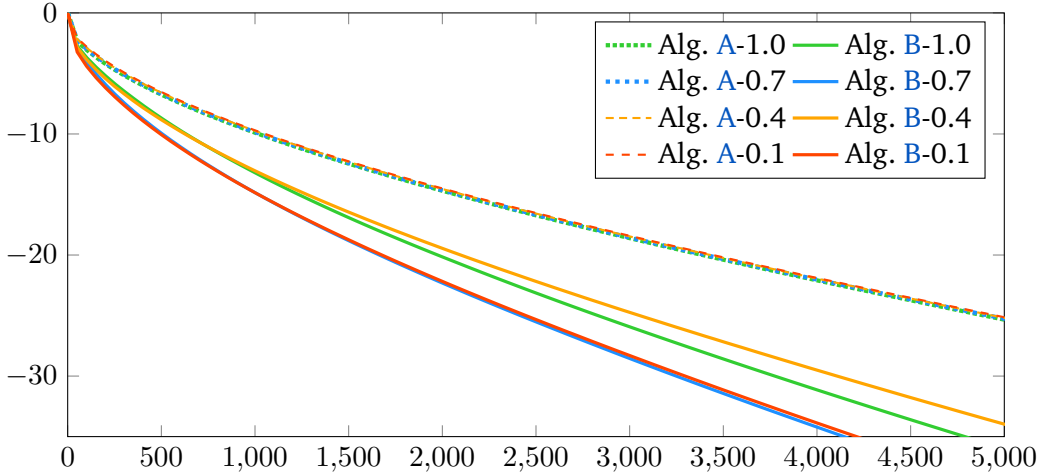 

\section{Numerical Experiments}
\label{sec:4}

We present two experiments which are reflective of our numerical
investigations in solving various problems using
Algorithms~\ref{algo:A} and \ref{algo:B}. The main objective is
to illustrate the block processing ability of the algorithms
(when implemented with full activation, i.e.,
$I_n=I_0$ and $K_n=K_0$, 
Algorithm~\ref{algo:B} was already shown in \cite{Siim19} to be
quite competitive compared to existing methods).

\subsection{Experiment 1: Group-Sparse Binary Classification}

We revisit the classification problem of \cite{Nume19}, which is
based on the latent group lasso formulation in machine learning
\cite{Jaco09}. Let $\{G_1,\ldots,G_m\}$ be a covering of
$\{1,\dots,d\}$ and define 
\begin{equation}
X=\menge{(x_1,\ldots,x_m)}{x_i\in\RR^d,\:
\operatorname{support}(x_i)\subset G_i}. 
\end{equation}
The sought vector is 
$\widetilde{y}=\sum_{i=1}^m\widetilde{x}_i$, where
$(\widetilde{x}_1,\ldots,\widetilde{x}_m)$ solves
\begin{equation}
\label{e:ml}
\minimize{(x_1,\ldots,x_m)\in X}{\sum_{i=1}^m\|x_i\|_{2}+
\sum_{k=1}^pg_k\Bigg(\sum_{i=1}^m\scal{x_i}{u_k}\Bigg)},
\end{equation}
with $u_k\in\RR^d$ and 
$g_k\colon\RR\to\RR\colon\xi\mapsto 10\max\{0,1-\beta_k\xi\}$,
where 
$\beta_k=\omega_k\operatorname{sign}(\scal{\overline{y}}{u_k})$
is the $k$th measurement of the true vector $\overline{y}\in\RR^d$
($d=10000$) and $\omega_k\in\{-1,1\}$ induces $25\%$ classification
error. There are $p=1000$ measurements and the
goal is to reconstruct the group-sparse vector $\overline{y}$. 
There are $m=1429$ groups. For every
$i\in\{1,\ldots,m-1\}$, each $G_i$ has $10$ consecutive integers
and an overlap with $G_{i+1}$ of length $3$. We obtain an
instance of \eqref{e:prob1}, where $\HH_i=\RR^{10}$,
$f_i=\|\cdot\|_2$, and 
$L_{k,i}=\scal{\mute}{\left.u_k\right|_{G_i}}$. The
auxiliary tasks for Algorithm~\ref{algo:A} (see
\ref{r:ABi}) are negligible \cite{Nume19}.
For each $\alpha\in\{0.1,0.4,0.7,1.0\}$, at iteration $n\in\NN$,
$I_n$ has $\lceil\alpha m\rceil$ elements and
the proximity operators of the scalar functions
$(g_k)_{1\leq k\leq p}$ are all used, i.e., $K_n=\{1,\ldots,p\}$.
We display in Fig.~\ref{fig:1} the normalized error versus the
epoch, that is, the cumulative number of activated blocks in
$\{1,\ldots,m\}$ divided by $m$.

\subsection{Experiment 2: Image Recovery}

\begin{figure}[t]
\centering
\begin{tabular}{@{}c@{}c@{}}
\includegraphics[width=6.2cm]{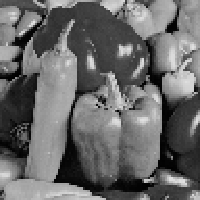}&
\hspace{0.2cm}
\includegraphics[width=6.2cm]{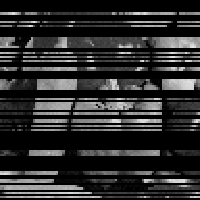}\\
\small{(a)} & \small{(b)}\\
\includegraphics[width=6.2cm]{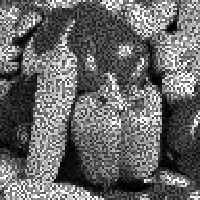}&
\hspace{0.2cm}
\includegraphics[width=6.2cm]{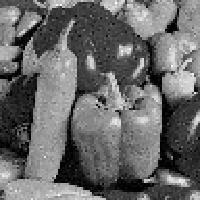}\\
\small{(c)} & \small{(d)}\\
\end{tabular} 
\caption{Experiment 2:
(a) Original $\overline{x}$.
(b) Observation $b$.
(c) Observation $c$.
(d) Recovery (all recoveries were visually 
indistinguishable).}
\label{fig:2}
\vskip -1mm
\end{figure}

We revisit the image interpolation problem of
\cite[Section~4.3]{Siim19}. The objective is to recover the image
$\overline{x}\in C=\left[0,255\right]^{N}$ ($N=96^2$) of
Fig.~\ref{fig:2}(a), given a noisy masked observation
$b=M\overline{x}+w_1\in\RR^N$ and a noisy blurred observation
$c=H\overline{x}+w_2\in\RR^N$. Here, $M$ masks all but $q=39$ rows
$(x^{(r_k)})_{1\leq k\leq q}$ of an image $x$, and $H$ 
is a nonstationary blurring operator, while $w_1$ and $w_2$
yield signal-to-noise ratios of $28.5$~dB and $27.8$~dB,
respectively. Since $H$ is sizable, we
split it into $s=384$ subblocks: for every 
$k\in\{1,\ldots,s\}$,
$H_{k}\in\RR^{24\times N}$ and the corresponding block of
$c$ is denoted $c_k$. The goal is to
\begin{equation}
\label{e:10}
\minimize{x\in C}{\|Dx\|_{1,2}+
10\sum_{k=1}^q\big\|x^{(r_k)}-b^{(r_k)}\big\|_2}
+5\sum_{k=1}^s\|H_kx-c_k\|_2^2,
\end{equation}
where $D\colon\RR^N\to\RR^N\times\RR^N$ models finite differences
and $\|\mute\|_{1,2}\colon
(y_1,y_2)\mapsto\sum_{j=1}^{N}\|(\eta_{1,j},\eta_{2,j})\|_2$.
Thus, \eqref{e:10} is an instance of Problem~\ref{prob:1}, where
$m=1$; $p=q+s+1$; for every $k\in\{1,\ldots,q\}$, 
$L_{k,1}\colon\RR^{N}\to\RR^{\sqrt{N}}\colon x\mapsto x^{(r_k)}$
and $g_k\colon y_k\mapsto 10\|y_k-b^{(r_k)}\|_2$;
for every $k\in\{q+1,\ldots,q+s\}$, $L_{k,1}=H_{k-q}$,
$g_k=5\|\mute-c_k\|_2^2$, and $g_p=\|\mute\|_{1,2}$; $L_{p,1}=D$;
$f_1\colon x\mapsto 0$ if $x\in C$; $\pinf$ if $x\not\in C$.
At iteration $n$, $K_n$ has $\lceil \alpha p\rceil$ elements,
where $\alpha\in\{0.1,0.4,0.7,1.0\}$.
The results are shown in Figs.~\ref{fig:2}--\ref{fig:3}, where
the epoch is the cumulative number of activated blocks in
$\{1,\ldots,p\}$ divided by $p$.

\subsection{Discussion}

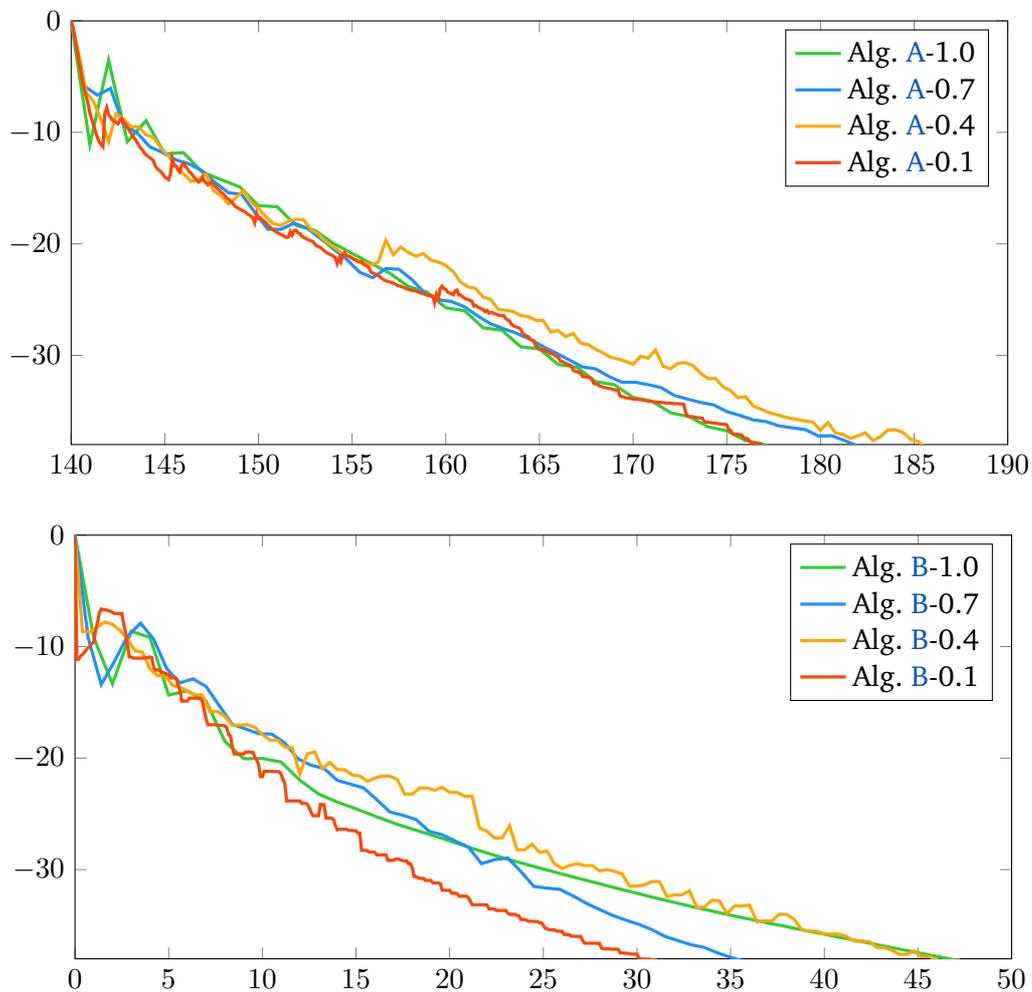
\begin{figure}
\centering
\begin{tikzpicture}[scale=1.00]
\begin{axis}[height=7.2cm,width=13.9cm, legend columns=1 
cell align={left}, xmin=140, xmax=190, ymin=-38, ymax=0.00]
\addplot [very thick,mark=none,color=pgreen]
table[x={epoch}, y={DR-err}] {figures2/ex2/DR-10.txt};
\addplot [very thick,mark=none,color=pblue]
table[x={epoch}, y={DR-err}] {figures2/ex2/DR-7.txt};
\addplot [very thick, mark=none,color=pyellow]
table[x={epoch}, y={DR-err}] {figures2/ex2/DR-4.txt};
\addplot [very thick, mark=none,color=pred]
table[x={epoch}, y={DR-err}] {figures2/ex2/DR-1.txt};
\legend{%
Alg.~\ref{algo:A}-1.0,
Alg.~\ref{algo:A}-0.7,
Alg.~\ref{algo:A}-0.4,
Alg.~\ref{algo:A}-0.1}
\end{axis}
\end{tikzpicture}\\[4mm]
\begin{tikzpicture}[scale=1.00]
\centering
\begin{axis}[height=7.2cm,width=13.9cm, legend columns=1 
cell align={left}, xmin=0, xmax=50, ymin=-38, ymax=0.00]
\addplot [very thick,mark=none,color=pgreen]
table[x={epoch}, y={CE-err}] {figures2/ex2/CE-10.txt};
\addplot [very thick,mark=none,color=pblue]
table[x={epoch}, y={CE-err}] {figures2/ex2/CE-7.txt};
\addplot [very thick, mark=none,color=pyellow]
table[x={epoch}, y={CE-err}] {figures2/ex2/CE-4.txt};
\addplot [very thick, mark=none,color=pred]
table[x={epoch}, y={CE-err}] {figures2/ex2/CE-1.txt};
\legend{%
Alg.~\ref{algo:B}-1.0,
Alg.~\ref{algo:B}-0.7,
Alg.~\ref{algo:B}-0.4,
Alg.~\ref{algo:B}-0.1}
\end{axis}
\end{tikzpicture}
\caption{Normalized error
$20\log_{10}(\|\boldsymbol{x}_n-\boldsymbol{x}_{\infty}\|
/\|\boldsymbol{x}_0-\boldsymbol{x}_{\infty}\|)$ (dB)
versus epoch count in Experiment 2. 
Top: Algorithm~\ref{algo:A}. The horizontal axis
starts at 140 epochs to account for the 
auxiliary tasks (see \ref{r:ABi}).
Bottom: Algorithm~\ref{algo:B}. The computational load
per epoch for Algorithm~\ref{algo:B} was about twice that of
Algorithm~\ref{algo:A}.}
\label{fig:3}
\vskip -3.5mm
\end{figure}

Our first finding is that, for both Algorithms~\ref{algo:A} and
\ref{algo:B}, even when full activation is computationally 
possible, it may not be the best strategy (see Figs.~\ref{fig:1}
and \ref{fig:3}). Second, \ref{r:ABi}--\ref{r:ABvii} and our
experiments suggest that \ref{algo:B} is preferable to
\ref{algo:A}. Let us add that, in general, \ref{algo:A}
does not scale as well as \ref{algo:B}. For instance, in
Experiment~2, if the image size scales up, \ref{algo:B}
can still operate since it involves only individual applications of
the local $L_{k,i}$ operators, while \ref{algo:A} becomes
unmanageable because of the size of the $Q_j$ operators (see
\ref{r:ABi} and \cite{Bric19}). 

\end{document}